\documentclass[a4paper,twoside,12pt]{article}
\usepackage[english]{babel}
\usepackage{graphicx}
\usepackage{bbding}
\usepackage[T1]{fontenc}
\usepackage{kpfonts}
\usepackage{fancyhdr}
\usepackage{amsmath,amssymb,color,epsfig}
\usepackage{mathrsfs}
\usepackage{eufrak}
\usepackage{dsfont}
\usepackage{upgreek}
\usepackage{fancyhdr}

\newtheorem{theorem}{Theorem}[section]
\newtheorem{proposition}[theorem]{Proposition}

\newtheorem{remark}[theorem]{Remark}
\newtheorem{lemma}[theorem]{Lemma}
\newtheorem{prediction}[theorem]{Prediction}

\def\1{\mathds{1}}

\title{Sobolev Freud polynomials}
\date{}
\begin{document}
\author{Mohamed BOUALI}
\maketitle
\begin{abstract} We investigate the uniform asymptotic of some Sobolev orthogonal polynomials. Three  term recurrence relation is given, moreover we give a recurrence relation between the so-called Sobolev orthogonal polynomials and Freud orthogonal polynomials.
\end{abstract}
\section{Introduction}
During the past few years, orthogonal polynomials with respect to an
inner product involving derivatives (so-called Sobolev orthogonal polynomials) have been the
object of increasing number of works (see, for instance [1], [5], [6], [4], [7], [8]). Recurrence
relations, asymptotics, algebraic, differentiation properties and zeros for various families of
polynomials have been studied. In this paper we study a connection between a particular case
of non-standard orthogonal polynomials.

For $\lambda_1,\lambda_2\geq 0$, we defined the inner product
$$\langle f,g\rangle_S=\int_{\Bbb R}f(x)g(x)e^{-x^4}dx+\lambda_1f(0)(x)g(0)+\lambda_2f'(0)(x)g'(0),$$
We denote also by $||\cdot||_S$ the norm associate to the inner product $\langle .,.\rangle_S$. Let $Q_n$ be the sequence of orthogonal polynomial with respect to $\langle .,.\rangle_S$. We denoted $\displaystyle \widehat k_n=||Q_n||^2_S=\langle Q_n, Q_n\rangle_s$ .

Let $P_n$ be the sequence of monic polynomials orthogonal with respect to the inner product $\displaystyle \langle f,g\rangle_F=\int_{-\infty}^{+\infty}f(x)g(x)e^{-x^4}dx$
:  They
have been considered by Nevai [14,15]. These polynomials satisfy a three-term
recurrence relation
$$xP_n(x)=P_{n+1}(x)+c_nP_{n-1}(x),$$
with initial conditions $P_0(x)=1$ and $P_1(x)=x$; where the parameters $c_n$ satisfy a
non-linear recurrence relation (see [4])
$$n=4c_n(c_{n+1}+c_n+c_{n-1}),\quad n\geq 1,$$
with $c_0=0$ and $c_1=\Gamma(3/4)/\Gamma(1/4)$. Moreover the polynomial $P_n$ satisfies the recurrence relation
$$P_n'(x)=nP_{n-1}(x)+d_nP_{n-3}(x),\quad n\geq 3$$
where $\displaystyle d_n=4{k_n}/{k_{n-3}},$ with
$\displaystyle k_n=||P_n||^2_F=\int_{-\infty}^{+\infty}\big(P_n(x)\big)^2e^{-x^4}dx.$
One can see from the three-term recurrence relation that \begin{equation}k_n=c_nk_{n-1},\end{equation}
and $$d_n=4c_nc_{n-1}c_{n-2}.$$
\begin{lemma} \
\begin{enumerate}
\item For all polynomials $P, Q$,
$\displaystyle \langle x^mP, Q\rangle_S=\langle P, x^mQ\rangle_S=\langle x^mP, Q\rangle_F=\langle P, x^mQ\rangle_F$,
where $m\geq 1$ if $\lambda_2=0$, and $m\geq 2$ if $\lambda_2>0$.
\item For a polynomial $Q$, we denote $\widetilde Q(x)=Q(-x)$. For all polynomials $P, Q$, we have $\langle\widetilde Q, P\rangle_S= \langle Q, \widetilde P\rangle_S$.
\item $Q_n(-x)=(-1)^nQ_n(x)$.
\end{enumerate}
\end{lemma}

{\bf Proof.}\\
The proof of the   first item is easy.\\
2) Using the symmetry of the Freud inner product $\langle\widetilde Q, P\rangle_F=\langle Q, \widetilde P\rangle_F$ and the fact that ${\widetilde P}'(0)=-P'(0)$
 $$\begin{aligned}\langle\widetilde Q, P\rangle_S&=\langle\widetilde Q, P\rangle_F+\lambda_1Q(0)P(0)-\lambda_2Q'(0)P'(0)\\&=\langle Q, \widetilde P\rangle_F+\lambda_1Q(0)P(0)-\lambda_2Q'(0)P'(0)\\
&=\langle Q, \widetilde P\rangle_F+\lambda_1Q(0){\widetilde P}(0)+\lambda_2Q'(0){\widetilde P}'(0)\\&=\langle Q, \widetilde P\rangle_S\end{aligned}.$$
3) From the first step we have  by orthogonality $\langle\widetilde Q_n, P\rangle_S=0$, for all polynomials $P$ with deg$(P)\leq n-1$.

Hence $\widetilde {Q_n}(x)=\alpha_nQ_n(x)$, equaling the leading coefficient we obtain $\alpha_n=(-1)^n$.

\section{Case $\lambda_2=0$}
\begin{proposition} The polynomials $P_n$ and $Q_n$ are related by
$$xP_n(x)=Q_{n+1}(x)+a_nQ_{n-1}(x),\quad n\geq 1$$
$$xQ_n(x)=P_{n+1}(x)+b_nP_{n-1}(x),\quad n\geq 1$$
$Q_0(x)=1$, $Q_1(x)=x$,
where $\displaystyle a_n=\frac{k_n}{\widehat k_{n-1}},$ $\displaystyle b_n=\frac{\widehat k_n}{k_{n-1}},$.
\end{proposition}
{\bf Proof.} Since, let  write \begin{equation}\label{1}xP_n(x)=Q_{n+1}(x)+\sum_{k=0}^{n}\alpha_kQ_k(x),\end{equation}
By orthogonality one gets
$$||Q_k||^2_S\alpha_k=\int_{-\infty}^{+\infty}P_n(x)Q_k(x)xe^{-x^4}dx,$$
Since for $k\leq n-2$, by orthogonality  the integral vanishes, moreover using the symmetry of the inner product, one as
$$P_n(-x)=(-1)^nP_n(x),\quad Q_k(-x)=(-1)^kQ_k(x),$$
hence $\alpha_n=0$, and $a_n=\alpha_{n-1}$.
$$\displaystyle\widehat k_{n-1} \alpha_{n-1}=\int_{-\infty}^{+\infty} P_n(x)xQ_{n-1}(x)e^{-x^4}dx=k_n,$$
where we used $xQ_{n-1}(x)=P_n(x)+...$
The second statement can be proved by a same argument.

\begin{proposition} \
\begin{enumerate}

\item $\displaystyle\frac{c_{n+2}c_{n+1}}{a_{n+2}}+a_n=c_{n+1}+c_n.$
\item $\displaystyle a_{n+1}b_n=c_{n+1}c_n$
\item $\displaystyle\lim_{n\to\infty}\frac{a_n}{\sqrt n}=\displaystyle\lim_{n\to\infty}\frac{b_n}{\sqrt n}=\frac{1}{2\sqrt 3}$.
\end{enumerate}
\end{proposition}
{\bf Proof.}\\
1) Since $xP_n(x)=Q_{n+1}+a_nQ_{n-1}(x),$ hence
$$||xP_n||^2_S=||Q_{n+1}||^2_S+a_n^2||Q_{n-1}||^2_S,$$
moreover $$||Q_{n-1}||^2_S=\frac{k_n}{a_n},$$
and $||x P_n||^2_S=||xP_n||^2_F$, from the tree-term recurrence  relation one gets
$$||xP_n||^2_F=||P_{n+1}||^2_F+c_n^2||P_{n-1}||^2_F=(c_{n+1}+c_n)k_n,$$
thus
$$\frac{k_{n+2}}{a_{n+2}}+a_{n}k_n=(c_{n+1}+c_n)k_n,$$
hence $$\frac{c_{n+2}c_{n+1}}{a_{n+2}}+a_{n}=c_{n+1}+c_n.$$
2) We saw that  $\displaystyle a_{n+1}=\frac{k_{n+1}}{\widehat k_{n}},$ $\displaystyle b_n=\frac{\widehat k_n}{k_{n-1}},$ and $k_n=c_n k_{n-1}$,hence
$$a_{n+1}b_n=\frac{k_{n+1}}{ k_{n-1}}=c_{n+1}c_n.$$
3) Let $\lambda_n=a_n\sqrt n$, $\sigma_n=\frac{c_{n+2}c_{n+1}}{n+1}$, $\delta_n=\frac{c_{n+1}+c_n}{\sqrt{n+1}}$.
Then we obtain, \begin{equation}\label{3}\frac{\sigma_n}{\lambda_{n+1}}+\frac{\sqrt{n-1}}{\sqrt{n+1}}\lambda_{n-1}=\delta_n,\end{equation}
Using the fact that $\displaystyle \lim_{n\to\infty}\frac{c_n}{\sqrt n}=\frac1{2\sqrt 3}$, (see for instance \cite{N}, \cite{b}) one gets, $$\lim_{n\to +\infty}\sigma_n=\frac1{12},\quad \lim_{n\to +\infty}\delta_n=\frac1{\sqrt 3}.$$
Moreover since $\lambda_n\geq 0$, and $$\lambda_n\leq\sqrt{\frac{n+2}{n}}\delta_{n+1},$$
thus the sequence  $\lambda_n$ is bounded. Let $\ell$ be the limit of a subsequence. It follows from equation (\ref{3}),
$$\frac{1}{12\ell}+\ell=\frac1{\sqrt 3 },$$
and the unique solution is $\ell=\frac1{2\sqrt 3}$. Hence the unique limit of a subsequence of $\lambda_n$ is $\ell=\frac1{2\sqrt 3}$, then the bounded sequence $\lambda_n$ converges to $\ell=\frac1{2\sqrt 3}$.

The same hold for $\frac{b_n}{\sqrt n}$ from the relation $a_{n+1}b_n=c_{n+1}c_n.$
\begin{theorem} The asymptotic behavior
$$\lim_{n\to\infty}\frac{Q_n(\sqrt[4]{n}x)}{P_n(\sqrt[4]{n}x)}=\sqrt[4]{12}\frac{x\varphi\big(\sqrt[4]{\frac34}x\big)}{1+\varphi^2\big(\sqrt[4]{\frac34}x\big)},$$
hold uniformly on compact subset of $\Bbb C\setminus[-\sqrt[4]{4/3},\sqrt[4]{4/3}]$, where\\ $\varphi(x)=x+\sqrt{x^2-1}$, with $\sqrt{x^2-1}>0$ for $x>1$, i.e., the conformal mapping of $\Bbb C\setminus[-1, 1]$ onto the exterior of
the closed unit disk.
\end{theorem}
{\bf Proof.} We saw that
$$xQ_n(x)=P_{n+1}(x)+b_nP_{n-1}(x).$$
It is well-known (see [16]) that from the three-term recurrence relation of non normalizing Freud polynomials
 $$xS_n(x)=\alpha_{n+1}S_{n+1}(x)+\alpha_nS_{n-1}(x),$$
 we can obtain asymptotic properties of the orthonormal polynomials $S_n$: Indeed, as
 $\displaystyle\lim_{n\to\infty}\frac{\alpha_n}{\sqrt[4]{n}}=\frac1{\sqrt[4]{12}}.$
 We deduce (see \cite{v})
 $$\lim_{n\to\infty}\frac{S_{n-1}\big(\sqrt[4]{n}x\big)}{S_n\big(\sqrt[4]{n}x\big)}=\frac{1}{\varphi\big(\sqrt[4]{\frac34}x\big)},$$
 uniformly on compact subsets of $\Bbb C\setminus[-\sqrt[4]{4/3},\;\sqrt[4]{4/3}]$. Then, for the monic Freud polynomial $P_n$ one gets
$$\lim_{n\to\infty}\sqrt[4]{n}\frac{P_{n-1}\big(\sqrt[4]{n}x\big)}{P_n\big(\sqrt[4]{n}x\big)}=\frac{\sqrt[4]{12}}{\varphi\big(\sqrt[4]{\frac34}x\big)},$$
uniformly on compact subsets of $\Bbb C\setminus[-\sqrt[4]{4/3},\;\sqrt[4]{4/3}]$.

Since $$\frac{Q_n(\sqrt[4]{n}x)}{P_n(\sqrt[4]{n}x)}=\frac1{x}\Big(\frac{P_{n+1}(\sqrt[4]{n}x)}{\sqrt[4]{n}P_n(\sqrt[4]{n}x)}+\frac{b_n}{\sqrt n}\frac{\sqrt[4]{n}P_{n-1}(x)}{P_n(x)}\Big).$$
As $n$ goes to infinity, one gets on every compact subsets of $\Bbb C\setminus[-\sqrt[4]{4/3},\;\sqrt[4]{4/3}]$,

$$\lim_{n\to\infty}\frac{Q_n(\sqrt[4]{n}x)}{{P_n\big(\sqrt[4]{n}x\big)}}=\frac1{x}\Big(\frac{\varphi(\sqrt[4]{\frac34}x)}{\sqrt[4]{12}}+
\frac{1}{\sqrt {12}}\frac{\sqrt[4]{12}}{\varphi(\sqrt[4]{\frac34})}\Big).$$
A simple computation  gives
$$\lim_{n\to\infty}\frac{Q_n(\sqrt[4]{n}x)}{{P_n\big(\sqrt[4]{n}x\big)}}=\frac1{\sqrt[4]{12}}\frac{1+\Big(\varphi\big(\sqrt[4]{\frac34}x\big)\Big)^2}
{x\varphi\big(\sqrt[4]{\frac34}x\big)}.$$
\begin{theorem}The polynomials $Q_n$ have all their zeros real and simple. For $n\geq 3$ the positive zeros of $Q_n$ interlace
with those of $P_n$.
\end{theorem}
{\bf Proof.} We distinguish two cases: the even and the odd one, respectively. The proofs are similar with slight differences.

{\bf Even case:} Let $x_{2m, k}$,$k=1,...,m$, be the positive zeros of $P_{2m}$ in increasing order, that is, $x_{2m, 1}<...<x_{2m, m}$.
First, we need to study the sign of the integrals
$$I_{2m, k}=\int_{-\infty}^{+\infty}Q_{2m}(x)\frac{P_{2m}(x)}{x^2-x^2_{2m, k}}e^{-x^4}dx,\quad m\geq 2,\; k=1,...,m,$$
We have
$$I_{2m, k}=\int_{-\infty}^{+\infty}Q_{2m}(x)\prod_{j=1, j\neq k}^{m-1}(x^2-x^2_{2m, j})e^{-x^4}dx=\sum_{r=0}^{m-1}b_r\int_{-\infty}^{+\infty}Q_{2m}(x)x^{2r}e^{-x^4}dx=b_0(k)\lambda_1Q_{2m}(0),$$
where $\displaystyle b_0(k)=(-1)^{m-1}\prod_{j=1, j\neq k}^{m-1}x^2_{2m,j}$,
and $\displaystyle Q_{2m}(0)=-a_{2m}Q_{2m-2}(0)=(-1)^m\prod_{k=2}^ma_{2k}$,\\($a_{2k}>0$), moreover
 ${\rm sign}(b_0)=(-1)^{m-1}$, hence
 \begin{equation}\label{17}{\rm sign}(I_{2m, k})=-1.\end{equation} On the other hand, using Gaussian quadrature in all the zeros of $P_{2m}$ and taking into account the symmetry of the polynomials $Q_{2m}$, the Christoffel numbers (see, for example, [7, p140])
 $$\mu_{2m, i}=\frac{1}{\sum_{j=0}^{2m-1}P^2_j(x_{2m, i})},\quad i=1,..,m, $$
 together with the fact
 $$\frac{P'_{2m}(x_{2m, k})}{2x_{2m,k}} =\prod_{j=1, j\neq k}^m\big(x^2_{2m, k}-x^2_{2m, j}\big),$$
 we get
 $$I_{2m, k}=\mu_{2m, k}Q_{2m}(x_{2m, k})\frac{P'_{2m, k}(x_{2m, k})}{2x_{2m, k}},$$
 and from $(\ref{17})$ we deduce
 $${\rm sign}(Q_{2m})=-{\rm sign}(P'_{2m, k}).$$
 Since $P'_{2m}(x)$ has opposite sign in two consecutive zeros of $P_{2m}(x)$, we deduce that it also occurs for $Q_{2m}(x)$, and therefore $Q_{2m}(x)$ has one zero in each interval $(x_{2m, k}, x_{2m+1, k}),\;k=1,...,m-1$ (and from the symmetry it has one zero in each interval $(-x_{2m+1, k}, -x_{2m, k}),\;k=1,...,m-1$. Thus $Q_{2m}(x)$ has at least $2m-2$ real and simple zeros interlacing with those of $P_{2m}(x)$. Finally, as $P'_{2m}(x_{2m, 2m})>0$ then $Q_{2m}(x_{2m, 2m})<0$ and since $Q_{2m}(x)$ is monic we deduce the existence
of one zero of $Q_{2m}(x)$ in $(x_{2m, m}, +\infty)$ and another zero in $( -\infty, -x_{2m, m})$, which  complete the result for the even case.\\
{\bf Odd case:} Let $m\ge 2$, $0<x_{2m+1,1}<...<x_{2m+1, 2m}$ be the positive simple zeros of $P_{2m+1}$, since $P_{2m+1}(0)= Q_{2m+1}(0)=0$, let define the integral
$$I_{2m+1, k}=\int_{-\infty}^{+\infty}Q_{2m+1}(x)\frac{P_{2m+1, k}(x)}{x^2(x^2-x^2_{2m+1, k})}e^{-x^4}dx$$
 hence
$$I_{2m+1, k}=\int_{-\infty}^{+\infty}\frac{Q_{2m+1}(x)}{x}\prod_{j=1, j\neq k}^{m}(x^2-x^2_{2m, j})e^{-x^4}dx,$$
$$I_{2m+1, k}=\sum_{r=0}^{m-1}b_r(k)\int_{-\infty}^{+\infty}\frac{Q_{2m+1}(x)}{x}x^{2r}e^{-x^4}dx,$$
Since for $ 1\leq r\leq m-1$, $$\int_{-\infty}^{+\infty}\frac{Q_{2m+1}(x)}{x}x^{2r}e^{-x^4}dx=\langle Q_{2m+1}, x^{2r-1}\rangle_S-\lambda_1(Q_{2m+1}x^{2r-1})\mid_{x=0},$$
hence by orthogonality one gets for $ 1\leq r\leq m-1$ $$\int_{-\infty}^{+\infty}\frac{Q_{2m+1}(x)}{x}x^{2r}e^{-x^4}dx=0.$$
Thus \begin{equation}\label{21}I_{2m+1, k}=b_0(k)\int_{-\infty}^{+\infty}\frac{Q_{2m+1}(x)}{x}e^{-x^4}dx,\end{equation}
since $$xP_{2m}(x)=Q_{2m+1}(x)+a_{2m}Q_{2m-1}(x),$$
hence $$\int_{-\infty}^{+\infty}P_{2m}(x)e^{-x^4}dx=\int_{-\infty}^{+\infty}\frac{Q_{2m+1}(x)}{x}e^{-x^4}dx+a_{2m}\int_{-\infty}^{+\infty}\frac{Q_{2m-1}(x)}{x}e^{-x^4}dx$$
thus by orthogonality $\int_{-\infty}^{+\infty}P_{2m}(x)e^{-x^4}dx=0$, and
$$\int_{-\infty}^{+\infty}\frac{Q_{2m+1}(x)}{x}e^{-x^4}dx=-a_{2m}\int_{-\infty}^{+\infty}\frac{Q_{2m-1}(x)}{x}e^{-x^4}dx$$
and
\begin{equation}\label{22}\int_{-\infty}^{+\infty}\frac{Q_{2m+1}(x)}{x}e^{-x^4}dx=2\Gamma(\frac54)(-1)^{m}\prod_{k=1}^ma_{2k}.\end{equation}
 Moreover \begin{equation}\label{23}b_0(k)=(-1)^{m-1}\prod_{j=1, j\neq k}^mx^2_{2m+1, j},\end{equation}
 from equations (\ref{21}), (\ref{22}) and (\ref{23}) one gets
  $${\rm sign}(I_{2m+1, k})=-1.$$
 The rest of the proof is as in the even case.
\section{Case $\lambda_2\neq 0$}
\begin{proposition}
For all $n\geq 1$,
$$xP_{2n-1}(x)=Q_{2n}(x)+a_nQ_{2n-2}(x),$$
$$x^2P_n(x)=Q_{n+2}(x)+b_nQ_{n}(x)+\alpha_nQ_{n-2}(x),$$
$$x^2Q_n(x)=P_{n+2}(x)+\sigma_nP_n(x)+\delta_n P_{n-2}(x).$$
with $Q_0(x)=1$, $Q_1(x)=x$, where,  $a_n=\frac{k_{2n-1}}{\widehat k_{2n-2}}$, $\alpha_n=\frac{k_n}{\widehat k_{n-2}}$, 
$b_n=\frac{\langle x^2P_ n,\,Q_n\rangle_S}{\langle Q_n,\,Q_n\rangle_S}$, 
 $\displaystyle \delta_n=\frac{\widehat k_n}{k_{n-2}}$, $\displaystyle\sigma_n=b_n\frac{\widehat k_n}{k_n}$.
\end{proposition}
{\bf Proof.} The proof is as in proposition 1.1.
\begin{proposition}
For all $n\geq 1$,
\begin{enumerate}
\item $\displaystyle\frac{c_{2n+1}c_{2n}}{a_{n+1}}+a_n=c_{2n}+c_{2n-1}.$
\item $\displaystyle c_{n+2}c_{n+1}+c_nc_{n-1}+(c_{n+1}+c_n)^2=\frac{c_{n+4}c_{n+3}c_{n+2}c_{n+1}}{\alpha_{n+4}}+b_n^2\frac{c_{n+2}c_{n+1}}{\alpha_{n+2}}+\alpha_n.$
    \item $\displaystyle \sigma_n=\frac{n}{4c_n}+c_{n-2}-b_{n-2}.$
    \item $\displaystyle c_{n+2}c_{n+1}\frac{\sigma_{n+2}}{b_{n+2}}+b_n\sigma_n+\alpha_n=\frac n2+\frac14.$
 \item $\displaystyle \sigma_n=\frac{\delta_n b_n }{c_nc_{n-1}}.$
    \item $\displaystyle\lim_{n\to\infty}\frac{a_n}{\sqrt {2n}}=\frac1{2\sqrt 3}$, $\displaystyle\lim_{n\to\infty}\frac{b_n}{\sqrt n}=\frac{1}{\sqrt 3}$, $\displaystyle\lim_{n\to\infty}\frac{\alpha_n}{n}=\frac{1}{12}$.
    \item $\displaystyle\lim_{n\to\infty}\frac{\sigma_n}{\sqrt n}=\frac1{\sqrt 3}$, $\displaystyle\lim_{n\to\infty}\frac{\delta_n}{n}=\frac1{12}$, $\displaystyle\lim_{n\to\infty}\frac{\widehat k_n}{k_n}=1$.

       \end{enumerate}
\end{proposition}
{\bf Proof.} \\
1) The first relation can be proved as the case $\lambda_2=0$.\\
2) From the second relation in the previous proposition and orthogonality one gets
$$||x^2P_n||^2_S=\widehat k_{n+2}+b_n^2\widehat k_n+\alpha_n^2\widehat k_{n-2}\\
.$$
Using the fact that $k_n=c_nk_{n-1}$, and $\widehat k_n=\frac{k_{n+2}}{\alpha_{n+2}}.$ One gets
\begin{equation}\label{9}
||x^2P_n||^2_S=\big(\frac{c_{n+4}c_{n+3}c_{n+2}c_{n+1}}{\alpha_{n+4}}+b_n^2\frac{c_{n+2}c_{n+1}}{\alpha_{n+2}}+\alpha_n\big)k_n
\end{equation}
Since from the three term recurrence relation $xP_n(x)=P_{n+1}(x)+c_nP_{n-1}(x)$, and orthogonality we have
$$||xP_n||^2_S=||xP_n||^2_F=k_{n+1}+c^2_nk_{n-1},$$
and \begin{equation}\label{8}\begin{aligned}&||x^2P_n||^2_S=||x^2P_{n}||^2_F=||x P_{n+1}+x c_{n}P_{n-1}||^2_F\\
&=||x P_{n+1}||_F +c^2_{n}||xP_{n-1}||^2_F+2c_{n}\langle xP_{n+1},\,xP_{n-1}\rangle_F\\
&=k_{n+2}+c_{n+1}^2k_{n}+c^2_{n}k_{n}+c_{n}^2c^2_{n-1}k_{n-2}+2c_{n}k_{n+1}\\
&=(c_{n+2}c_{n+1}+c_{n+1}^2+c^2_{n}+c_{n}c_{n-1}+2c_{n}c_{n+1})k_n\\
&=\big(c_{n+2}c_{n+1}+(c_{n+1}+c_{n})^2+c_{n}c_{n-1}\big)k_{n}.
\end{aligned}\end{equation}
From equation (\ref{9}), (\ref{8}) one gets the desired result.

3) By orthogonality one gets
$$\sigma_nk_n=\langle x^2Q_n,\, P_n\rangle_F,$$
since by definition of the Sobolev inner product we have
$$\begin{aligned}\langle x^2Q_n,\, P_n\rangle_F&=\langle x^2Q_n,\, P_n\rangle_S=\langle Q_n,\, x^2P_n\rangle_S\\
&=\langle x^2P_{n-2}-b_{n-2}Q_{n-2}-\alpha_{n-2}Q_{n-4},\, x^2P_n\rangle_S\\
&=\langle x^2P_{n-2},\, x^2P_n\rangle_S-b_{n-2}\langle Q_{n-2},\, x^2P_n\rangle_S-\alpha_{n-2}\langle Q_{n-4},\, x^2P_n\rangle_S
\end{aligned},$$
since $$\begin{aligned}\langle x^2P_{n-2},\, x^2P_n\rangle_S&=\int_{-\infty}^{+\infty}x^4P_{n-2}(x)P_n(x)e^{-x^4}dx\\
&=\frac14\int_{-\infty}^{+\infty}P_{n-2}(x)P_n(x)e^{-x^4}dx+\frac14\int_{-\infty}^{+\infty}xP'_{n-2}(x)P_n(x)e^{-x^4}dx\\
&+\frac14\int_{-\infty}^{+\infty}xP_{n-2}(x)P'_n(x)e^{-x^4}dx,
\end{aligned}$$
by orthogonality the first and the second integral vanished. Moreover
$$\begin{aligned}\int_{-\infty}^{+\infty}xP_{n-2}(x)P'_n(x)e^{-x^4}dx=&
\int_{-\infty}^{+\infty}xP_{n-2}(x)\big(nP_{n-1}(x)+d_nP_{n-3}(x)\big)e^{-x^4}dx\\
&=nk_{n-1}+d_nk_{n-2}.
\end{aligned}$$
$$\langle Q_{n-2},\, x^2P_n\rangle_S=\langle x^2Q_{n-2},\, P_n\rangle_F=k_n.$$
$$\langle Q_{n-4},\, x^2P_n\rangle_S=\langle x^2Q_{n-4},\, P_n\rangle_F=0.$$
Hence
$$\sigma_nk_n=\frac n4k_{n-1}+\frac14 d_nk_{n-2}-b_{n-2}k_n,$$
using the fact that $k_n=c_nk_{n-1}$, and $d_n=4c_nc_{n-1}c_{n-2}.$ Thus
  $$\sigma_n=\frac n{4c_n}+c_{n-2}-b_{n-2},$$

which complete the proof of the assertion.

4) From the second recurrence relation of the proposition we have
 $$\langle x^2P_n,\,x^2P_n\rangle_S=\widehat k_{n+2}+b^2_n\widehat k_n+\alpha^2_n\widehat k_{n-2}.$$
Moreover $\langle x^2P_n,\,x^2P_n\rangle_S=\langle x^2P_n,\,x^2P_n\rangle_F$, and
$$\begin{aligned}\langle x^2P_n,\,x^2P_n\rangle_F&=\int_{-\infty}^{+\infty}x^4P^2_n(x)e^{-x^4}dx\\
&=\frac14\int_{-\infty}^{+\infty}P^2_n(x)e^{-x^4}dx+\frac12\int_{-\infty}^{+\infty}xP'_n(x)P_n(x)e^{-x^4}dx\\
&=\frac14 k_n+\frac12nk_n,
\end{aligned}$$
thus, $$\frac{\widehat k_{n+2}}{k_n}+b^2_n\frac{\widehat k_n}{k_n}+\alpha^2_n\frac{\widehat k_{n-2}}{k_n}=\frac14+\frac12n,$$
since $\displaystyle \frac{\widehat k_n}{k_n}=\frac{\sigma_n}{b_n}$, $k_n=\alpha_{n}{\widehat k_{n-2}}$, and $k_n=c_nk_{n-1}$, it follows that
$$c_{n+2}c_{n+1}\frac{\sigma_{n+2}}{\alpha_{n+2}}+b_n\sigma_n+\alpha_n=\frac n2+\frac14.$$
Which give  the desired result.

5) It is easy deduced from the three relations $\displaystyle\delta_n=\frac{\widehat k_n}{k_{n-2}}$,  $\displaystyle\sigma_n=b_n\frac{\widehat k_n}{k_n}$, and $\displaystyle k_n=c_n k_{n-1}$.

6) The first limit can be proven in the same way of the case $\lambda_2=0$, in fact we saw
$$\frac{c_{2n+1}c_{2n}}{a_{n+1}}+a_n=c_{2n}+c_{2n-1},$$
using the same argument as in the case $\lambda_2=0$, one gets $\displaystyle\lim_{n\to+\infty}\frac{a_n}{\sqrt {2n}}=\ell=\frac1{2\sqrt 3}.$

 To obtain the second limit, one can see from the first equation that
 $$0\leq \frac{\alpha_n}{n}\leq \frac1n(c_{n+2}c_{n+1}+c_nc_{n-1}+(c_{n+1}+c_n)^2),$$
 since  the right hand side converge, hence $\frac{\alpha_n}{n}$ is bounded. And
 $$\big(\frac{b_n}{\sqrt n}\big)^2\leq \frac{c_{n+2}c_{n+1}+c_nc_{n-1}+(c_{n+1}+c_n)^2}{c_{n+2}c_{n+1}}\frac{\alpha_{n+2}}{n},$$
 moreover the sequence $\displaystyle \frac{c_n}{\sqrt n}$ converge and $\displaystyle\frac{\alpha_{n}}{n}$ is bounded hence the sequence $\displaystyle\frac{b_n}{\sqrt n}$ is bounded. Let $x$ be a limit of any subsequence of $\displaystyle\frac{b_{n_k}}{\sqrt {n_k}}$ and $\ell=\frac1{2\sqrt 3}$ the limit of the sequence $\displaystyle \frac{c_n}{\sqrt n}$. Since from the third and fourth items we have
 $$\sigma_n=\frac{n}{4c_n}+c_{n-2}-b_{n-2},$$
 and $$\alpha_n=\frac n2+\frac14-c_{n+2}c_{n+1}\frac{\sigma_{n+2}}{b_{n+2}}-b_n\sigma_n.$$
 Thus $$\alpha_n=\frac n2+\frac14-c_{n+2}c_{n+1}\frac1{b_{n+2}}(\frac{n+2}{4c_{n+2}}+c_{n}-b_{n})-b_n(\frac{n}{4c_n}+c_{n-2}-b_{n-2}).$$
 Substitute the expression of $\alpha_n$  in the first item of the proposition and letting $k$ to infinity and use the fact that  $\displaystyle\frac{b_{n_k}}{\sqrt {n_k}}$ converge to $x$ and  $\displaystyle \frac{c_{n_k}}{\sqrt {n_k}}$ converge to $\ell$, one gets
 $$\begin{aligned}-6\ell^2&=\frac{\ell^4}{-\frac12+x(\frac1{4\ell}+\ell-x)+\frac1{x}(\frac\ell 4+\ell^2-\ell x)}+\frac{x^2\ell^2}{-\frac12+x(\frac1{4\ell}+\ell-x)+\frac1{x}(\frac\ell 4+\ell^2-\ell x)}\\&+\frac12+x(\frac1{4\ell}+\ell-x)+\frac1{x}(\frac\ell 4+\ell^2-\ell x),\end{aligned}$$
 substitute the value of $\ell=\frac1{2\sqrt3}$ we obtain the following equation
 \begin{equation}\label{s}\frac{\left(1+12 x^2\right) \left(1-4 \sqrt{3} x+18 x^2-12 \sqrt{3} x^3+9 x^4\right)}{3 x \left(2 \sqrt{3}-21 x+24 \sqrt{3} x^2-36 x^3\right)}=\frac{\frac19\left(x-\frac1{\sqrt 3}\right)^4\left(1+12x^2\right)}{3 x \left(2 \sqrt{3}-21 x+24 \sqrt{3} x^2-36 x^3\right)}=0,\end{equation}
 Now we prove that such equation is correctly defined. Since
 the roots of the polynomial  $ x (2 \sqrt{3}-21 x+24 \sqrt{3} x^2-36 x^3)$ are $0,\; \ell=\frac1{2\sqrt 3}$ and two complex roots, moreover $\frac{b_{n_k}}{\sqrt n_{k}}$ is a real sequence.

 {\bf First case} $x=\ell=\frac1{2\sqrt 3}$. Since,
 $$\sigma_{n_k}=\frac{n_k}{4c_{n_k}}+c_{n_k-2}-b_{n_k-2},$$
 hence $\displaystyle\frac{\sigma_{n_k}}{n_k}$ converge to $\displaystyle\frac1{4\ell}=\frac{\sqrt 3}{2}.$ Moreover
 $\displaystyle\frac{\alpha_n}{\sqrt n}$ is bounded, then we can subtracted from $\displaystyle\frac{\alpha_n}{\sqrt n}$ a sequences which converges to some $y$, using the convergence of the sequence $\displaystyle\frac{b_{n_k}}{\sqrt n_k}$ to $x$, and equations 2) and 4) of the proposition,
 one gets from equation 2)
 $$6\ell^2=\frac{\ell^4}{y}+\frac{\ell^4}{y}+y,$$
 and from 4)
 $$y=0,$$
 which give a contradiction.

  {\bf Second case}. $x=0$. In such a case, we obtain from statements 3) and 4) and boundedness of $\frac{\alpha_n}{\sqrt n}$, one gets
  $$\lim_{k\to\infty}\frac{\sigma_{n_k}}{n_k}=\frac1{4\ell}+\ell,$$ and
 $$\lim_{k\to\infty}\frac{\sigma_{n_k}}{n_k}=0,$$
This give a contradiction.

 The only real solution of equation (\ref{s}) is $\displaystyle x=\frac1{\sqrt 3}$. Hence the unique accumulation point of the bounded sequence $\displaystyle \frac{b_n}{\sqrt n}$ is  $\displaystyle\frac1{\sqrt 3}$.

{\bf The sequence $\sigma_n$.} Since $\displaystyle \sigma_n=\frac{n}{4c_n}+c_{n-2}-b_{n-2}$, hence
$$\lim_{n\to\infty}\frac{\sigma_n}{\sqrt n}=\frac1{4\ell}+\ell-2\ell=\frac1{\sqrt 3}.$$
{\bf The sequence $\alpha_n$}. Letting $n$ to infinity in the following equation
$$\displaystyle c_{n+2}c_{n+1}\frac{\sigma_{n+2}}{b_{n+2}}+b_n\sigma_n+\alpha_n=\frac n2+\frac14,$$
one gets
$$\lim_{n\to\infty}\frac{\alpha_n}{n}=\frac12-\ell^2-4\ell^2=\frac1{12}.$$
7)
{\bf The sequence $\delta_n$.} From the relation
$$\displaystyle \delta_n=\frac{c_nc_{n-1}\sigma_n}{b_n},$$
as $n$ goes to infinity one gets
$$\lim_{n\to\infty}\frac{\delta_n}{n}=\ell^2=\frac1{12}.$$
{\bf The sequence $\frac{\widehat k_n}{k_n}$.} Since
$$\lim_{n\to\infty}\frac{\widehat k_n}{k_n}=\lim_{n\to\infty}\frac{\sigma_n}{b_n}=1.$$
Which complete the proof.
\begin{theorem}
The asymptotic behavior
$$\lim_{n\to\infty}\frac{P_{n}(\sqrt[4]{n}x)}{Q_{n}(\sqrt[4]{n}x)}=
2\sqrt{3}\Big(\frac{x\varphi\big(\sqrt[4]{\frac34}x\big)}{1+\varphi^2\big(\sqrt[4]{\frac34}x\big)}\Big)^2,$$
hold uniformly on compact subset of $\Bbb C\setminus[-\sqrt[4]{4/3},\sqrt[4]{4/3}]$, where\\ $\varphi(x)=x+\sqrt{x^2-1}$, with $\sqrt{x^2-1}>0$ for $x>1$, i.e., the conformal mapping of $\Bbb C\setminus[-1, 1]$ onto the exterior of
the closed unit disk.
\end{theorem}
{\bf Proof.} Sine from the third relation in proposition 2.2, we have
$$x^2Q_n(x)=P_{n+2}(x)+\sigma_nP_n(x)+\delta_n P_{n-2}(x),$$
thus
$$x^2\frac{Q_n(\sqrt[4]{n}x)}{P_n(\sqrt[4]{n}x)}=\frac{P_{n+2}(\sqrt[4]{n}x)}{\sqrt[4]{n}P_{n+1}(\sqrt[4]{n}x)}
\frac{P_{n+1}(\sqrt[4]{n}x)}{\sqrt[4]{n}P_{n}(\sqrt[4]{n}x)}
+\frac{\sigma_n}{\sqrt n}
+\frac{\delta_n}{n} \frac{\sqrt[4]{n}P_{n-2}(\sqrt[4]{n}x)}{P_{n-1}(\sqrt[4]{n}x)}\frac{\sqrt[4]{n}P_{n-1}(\sqrt[4]{n}x)}{P_{n}(\sqrt[4]{n}x)},$$
using the fact that $\displaystyle\frac{\sigma_n}{\sqrt n}\to \frac1{\sqrt 3}$,  $\displaystyle\frac{\delta_n}{\sqrt n}\to \frac1{12}$ and
$$\lim_{n\to\infty}\frac{\sqrt[4]{n}P_{n-1}(\sqrt[4]{n}x)}{P_{n}(\sqrt[4]{n}x)}=\frac{\sqrt[4]{12}}{\varphi\big(\sqrt[4]{\frac34}x\big)},$$
uniformly on compact subset of $\Bbb C\setminus[-\sqrt[4]{4/3},\sqrt[4]{4/3}]$, one gets
$$\lim_{n\to\infty}x^2\frac{Q_n(\sqrt[4]{n}x)}{P_n(\sqrt[4]{n}x)}=\frac{\Big(\varphi\big(\sqrt[4]{\frac34}x\big)\Big)^2}{\sqrt{12}}+\frac1{\sqrt  3}+\frac{1}{\sqrt 12\Big(\varphi\big(\sqrt[4]{\frac34}x\big)\Big)^2}=\frac{\Big(1+\Big(\varphi\big(\sqrt[4]{\frac34}x\big)\Big)^2\Big)^2}{\sqrt 12\Big(\varphi\big(\sqrt[4]{\frac34}x\big)\Big)^2},$$
which complete  the proof.
\begin{remark} \
\begin{enumerate}
\item For all $n$ the zeros of $Q_{2n}$ interlaces with the zeros of $P_{2n}$.
\item What can say about the zeros of the polynomial $Q_{2n+1}$ compared  with those of $P_{2n+1}$?.
\end{enumerate}
\end{remark}
The proof of the first item is like the proof of the case $\lambda_2=0$.

Consider the non monic Sobolev orthogonal polynomials $\widehat Q_n$, with norm equal one.
$\displaystyle\widehat Q_n(x)=c_n x^n+b_{n-2}x^{n-2}+....$.
\begin{proposition} \
\begin{enumerate}
\item For all $n\geq 0$, the polynomial $\widehat Q_n$, satisfies the three term recurrence relation
$$\displaystyle x^2\widehat Q_n(x)=\alpha_n\widehat Q_{n+2}+\beta_n\widehat Q_n(x)+\alpha_{n-2}\widehat Q_{n-2}(x),$$
with $Q_{-2}(x)=0$, $Q_{-1}(x)=0$.
\item The zeros of $\widehat Q_n$ interlaces with the zeros of $\widehat Q_{n-2}$.
\end{enumerate}
\end{proposition}
The proof of the proposition is like the proof for the classical orthogonal polynomial case.
\section{General case}
Let $r\in\Bbb N$, $\lambda_k\geq 0$ for $k\in\{1,...,r\}$, and defined the Sobolev inner product by
$$\langle f,\,g\rangle_S=\int_{-\infty}^{+\infty}f(x)g(x)e^{-x^4}dx+\sum_{k=0}^r\lambda_kf^{(k)}(0)g^{(k)}(0).$$
Let $Q_{n,r}$ be the monic Sobolev orthogonal polynomials with respect to the inner product defined above.
\begin{prediction} On every compact subset of $\Bbb C\setminus[-\sqrt[4]{4/3},\sqrt[4]{4/3}]$, we have
$$\lim_{n\to\infty}\frac{P_{n}(\sqrt[4]{n}x)}{Q_{n, r}(\sqrt[4]{n}x)}=
\Big(\sqrt[4]{12}\frac{x\varphi\big(\sqrt[4]{\frac34}x\big)}{1+\varphi^2\big(\sqrt[4]{\frac34}x\big)}\Big)^{r+1},$$
the convergence hold uniformly. 
\end{prediction}

Address:  College of Applied Sciences
   Umm Al-Qura University
  P.O Box  (715), Makkah,
  Saudi Arabia.\\
E-mail: bouali25@laposte.net \& mabouali@uqu.edu.sa
\end{document}